\documentclass[10pt,a4paper,oneside,fleqn,reqno]{amsart}

\usepackage[bookmarks=false]{hyperref}
\hypersetup{
	colorlinks=true,
 	linkcolor=blue,
	citecolor=blue,
	pdftitle={Bureau's Classification of Quadratic Differential Equations},
	pdfauthor={Adolfo Guillot},
}

\usepackage{ccicons}


\thanks{Published as: \textsc{A. Guillot}, On Bureau's Classification of Quadratic Differential Equations in Two Variables Free of Movable Critical points, \emph{Stud. Appl. Math.} \textbf{155} no. 2 (2025), 70083, \href{https://doi.org/10.1111/sapm.70083}{doi:10.1111/sapm.70083}}

\thanks{\ccby\, CC BY 4.0. This work is licensed under a \href{https://creativecommons.org/licenses/by/4.0/deed.en}{Creative Commons Attribution 4.0 License}}


\newtheorem{theorem}{Theorem}
\newtheorem{proposition}[theorem]{Proposition}
\newtheorem{corollary}[theorem]{Corollary}

\theoremstyle{remark}
\newtheorem{remark}[theorem]{Remark}

\newcommand{\PP}{\mathbf{P}}
\newcommand{\ZZ}{\mathbf{Z}}
\newcommand{\CC}{\mathbf{C}}
\newcommand{\dd}{\mathrm{d}}

\title[Quadratic differential equations]{On Bureau's classification of quadratic differential equations in two variables free of movable critical points}

\author{Adolfo Guillot}\address{Instituto de Matem\'{a}ticas, Universidad Nacional Aut\'{o}noma de M\'{e}xico, Ciudad Universitaria  04510, Mexico City,  Mexico}

\email{adolfo.guillot@im.unam.mx}

\begin{document}
\begin{abstract} As part of the efforts aimed at extending Painlev\'e and Gambier's work on second-order equations in one variable to first-order ones in two, in 1981, Bureau classified the systems of ordinary quadratic differential equations in two variables which are free of movable critical points (which have the Painlev\'e Property). We revisit this classification, which we  complete by adding some cases overlooked by Bureau, and by correcting some of his arguments. We also simplify the canonical forms of some systems, bring the natural symmetries of others into their study, and investigate the birational equivalence  among some of the systems in the class. Lastly, we study the birational geometry of Okamoto's space  of initial conditions for Bureau's system VIII, in order to establish the sufficiency of some necessary conditions for the absence of movable critical points.	\end{abstract}

\maketitle

\section{Introduction} Around the beginning of the  twentieth century, Painlev\'e considered the problem of determining   the algebraic ordinary  differential equations in the complex domain whose general solution is univalent,  starting with those of first order, then those of the second one, and so on (see \cite[1\(^\text{re}\) le\c con]{painleve-lecons}, \cite{painleve-bsmf}). At the time, the problem was completely understood for first-order equations, a class that contains the equations solved by elliptic functions; their rich theory and widely applicable nature motivated the search for higher-order analogs.

For a complex-valued function of one complex variable, a direct source of multivaluedness is   the existence of \emph{critical points}, singular points around which at least two of its branches are permuted (see  \cite{painleve-acta}). For the solutions of differential equations, these points may be separated into \emph{fixed} (occurring at fixed values of the independent variable for all solutions of a given equation) and \emph{movable} (appearing at values of the independent variable depending upon the particular solution of the equation). In order to tackle the previous  problem, Painlev\'e proposed the study of a wider class of differential equations, those  whose solutions are free of movable critical points (``\'equations \emph{\`a points critiques fixes}''). This property acquired a relevance of its own, and is sometimes referred to as the \emph{Painlev\'e Property}. He saw this class as a natural extension of that of linear differential equations, and   developed tools for its study. His methods evidenced the utility of considering equations whose dependence upon the independent variable is only analytic (and not necessarily algebraic).

Painlev\'e and Gambier's  efforts  lead to a classification of the second-order equations of the form \(y''=R(y',y,t)\), with \(R\) rational in \(y\) and \(y'\) and analytic in the independent variable \(t\), which are free of movable critical points: up to fractional linear transformations in \(y\), depending analytically on \(t\), and reparametrizations of the independent variable, they belong to some fifty explicit classes of equations that are either integrated in an elementary way, or are one of six model equations, the \emph{Painlev\'e equations}, labeled I trough VI, which cannot be reduced to either linear or first-order equations  
\cite{painleve-acta}, \cite{gambier-acta} (see also~\cite[Ch.~XIV]{ince}, \cite{bureau-m1}, \cite[Ch. 12.2]{hille}). Painlev\'e anticipated the applicability of these functions, and his insight has been proven amply correct, as they  ``have arisen in a variety of important physical applications including statistical mechanics, random matrices, plasma physics, nonlinear waves, quantum gravity, quantum field theory, general relativity, nonlinear optics and fibre optics'' \cite{Clarkson-Painleve} (see also \cite[\S 7.1.6]{Ablowitz-Clarkson}, \cite[Ch.~10]{gromak_et_al}).

The results of Painlev\'e and Gambier have been generalized to higher-order and higher-degree equations in one variable; we refer to  \cite[App.~A]{painleve_handbook} for a presentation of some  of these generalizations.

One may also consider the problem of classifying systems of differential equations of the form
\begin{equation} \label{eq:int} 
\begin{split}
		y' & =P(y,z,t), \\ z' & =Q(y,z,t),
	\end{split} 
\end{equation} 
with \(P\) and \(Q\) rational in \(y\) and \(z\), and analytic in the independent variable \(t\), which are free of movable critical points. The problem solved by Painlev\'e and Gambier is a particular case of it. In the above problem, it is natural to consider equivalence  up to changes of coordinates which are analytic in the independent variable and birational in the dependent ones, as these changes of coordinates preserve both the form of the equations and the property of being free of movable critical points. Garnier, who attributes the problem to Goursat, wondered if such general systems could define transcendents (functions) beyond those that appear in Painlev\'e and Gambier's classification, and estimated this to be a difficult problem \cite{garnier}. One may restrict it by considering systems of the form (\ref{eq:int}) with \(P\) and \(Q\) polynomial in \(y\) and \(z\),  up to either affine or birational  equivalence (the first preserves the class of systems of a given degree; the second, the systems whose solutions are free of movable critical points). These polynomial equations already form a rich class, which naturally contains Painlev\'e's first two, and which, as shown by Malmquist \cite[\S 33]{malmquist-17}, contains the other Painlev\'e equations as well (see also~\cite{okamoto}). Goffar-Lombet investigated the question for polynomial systems having some particular homogeneous polynomials, of degrees four and five, as highest-degree terms, obtaining precise and complete results for them~\cite{goffar-lombet}. In 1980, Kimura and Matuda announced a classification of polynomial systems of the form~(\ref{eq:int}) free of movable critical points, of degree up to five~\cite{kimura-matuda}, which does not seem to have had a follow-up from their part. Shortly afterwards, in a series of articles published in the \emph{Bulletin de la Classe des Sciences} of the \emph{Acad\'emie Royale de Belgique} (\cite{bureau-8}, \cite{bureau-9}, \cite{bureau-10}), Bureau classified quadratic systems up to affine transformations. His results are summarized in Table III in~\cite[\S 28]{bureau-10} (he uses the term \emph{stable} to refer to systems that are free of movable critical points).  The aim of this article is to revisit and revise Bureau's classification:

\begin{theorem}\label{thm:main} Consider a system of differential equations of the form~(\ref{eq:int}), with \(P\) and \(Q\) polynomials in \(y\) and \(z\) of degree at most two (with at least one of them non-linear), whose coefficients are holomorphic functions defined in a domain \(U\subset \CC\). If the system is free of movable critical points then, for every \(t_0\) in \(U\), with the possible exception of a discrete set, there exists a neighborhood \(V\) of \(t_0\) within \(U\), and a transformation \(\phi:\CC^2\times V\to \CC^2\times \CC\),
\[(\zeta,t)\mapsto (\xi(\zeta,t),\tau(t))\]
with \(\xi\) affine in \(\zeta\) and holomorphic in \(t\), and \(\tau\) a biholomorphism onto its image, such that the system solved by \(\xi\), as a function of~\(\tau\), is the restriction to \(\CC^2\times \phi(V)\) of  one of those appearing in Table~\ref{table:main}.

Reciprocally, all systems in Table~\ref{table:main} are free of movable critical points. They are integrated either in an elementary way or via some Painlev\'e transcendent (also indicated in Table~\ref{table:main}).\end{theorem}

\begin{table} 
\thispagestyle{empty}
\begin{tabular}{clc}
System	 & Equations & Integration  \\	 \hline
\hyperref[sec:1]{\(\text{I}\)} &
\begin{tabular}{l} 
\(y'=-y^2+az\) \\ 
\(z'=-yz\)
\end{tabular} & ---
  \\ \hline
\hyperref[sec:2]{\(\text{II}\)} &
\begin{tabular}{l} 
	\(y'=y(y-2z)+a\) \\ 
	\(z'=-z^2\)
\end{tabular} & ---
  \\ \hline
\hyperref[sec:IX_0]{\(\text{III}\)} &
\begin{tabular}{l} 
	\(y'=y^2\) \\ 
	\(z'=2yz+ay\)
\end{tabular} &---
  \\ \hline
\hyperref[sec:IV]{\(\text{IV}\)} &
\begin{tabular}{l} 
	\(y'=-y^2+az\) \\ 
	\(z'=0\)
\end{tabular} & ---
  \\ \hline 
\hyperref[sec:V]{\(\text{V}_0\)} & 
\begin{tabular}{l} 
	\(y'=0\) \\ 
	\(z'=6y^2\)
\end{tabular} & ---
  \\ \hline
\hyperref[sec:V]{\(\text{V}\)} &
\begin{tabular}{l} 
	\(y'=z\) \\ 
	\(z'=6y^2+f\) \\
	\quad with \(f''=0\)
\end{tabular}  & Painlev\'e I
 \\ \hline
\hyperref[sec:6]{\(\text{VI}\)} &
\begin{tabular}{l} 
	\(y'=a\) \\ 
	\(z'=z(z+y)+b\)
\end{tabular} & ---
  \\ \hline
\hyperref[sec:7]{\(\text{VII}_0\)} &
\begin{tabular}{l} 
	\(y'=0\) \\ 
	\(z'=yz+a\)
\end{tabular} & ---
  \\ \hline
\hyperref[sec:7]{\(\text{VII}\)} &
\begin{tabular}{l} 
	\(y'=z\) \\ 
	\(z'=yz+a\)
\end{tabular} & ---
  \\ \hline
\begin{tabular}{c}\hyperref[sec:VIII]{\(\text{VIII}(n)\)} \\ (\(n\geq 1\)) \end{tabular} &
\begin{tabular}{l} 
	\(y'=y[y-(n+2)z]+a\) \\ 
	\(z'=-z(ny+z)+b \) \\
\quad 	with a condition of order \(n-1\) (Section \ref{sec:VIII})
\end{tabular} & ---
  \\ \hline
\begin{tabular}{c} \hyperref[sec:IX_0]{\(\text{IX.A}_0(n)\)} \\ (\(n>1\)) \end{tabular} &
\begin{tabular}{l} 
	\(y'=-y^2\) \\ 
	\(z'=-(n+1)yz+ay\)
\end{tabular} & ---
  \\ \hline 
\hyperref[sec:IX.-3]{\(\text{IX.A}(3)\)} &
\begin{tabular}{l} 
	\(y'=-y^2+z\) \\ 
	\(z'=-4yz+a\)
\end{tabular}  & ---
  \\ \hline
\begin{tabular}{c} \hyperref[sec:IX_0]{\(\text{IX.B}_0(n)\)} \\ (\(n>1\)) \end{tabular} &
\begin{tabular}{l} 
	\(y'=-y^2\) \\ 
	\(z'=(n-1)yz+py\) \\
	\quad with \(\dd^{n-1}p/\dd t^{n-1}=0\)
\end{tabular} & ---
  \\ \hline
\hyperref[sec:IX.2]{\(\text{IX.B}(2)\)} &
\begin{tabular}{l} 
	\(y'=-y^2+z+12q\) \\ 
	\(z'=yz\) \\
		\quad with \((q''-6q^2)''=0\)
\end{tabular} & Painlev\'e I
  \\ \hline
\hyperref[sec:IX(3)]{\(\text{IX.B}(3)\)} &
\begin{tabular}{l} 
	\(y'=-y^2+z-\frac{1}{2}f\) \\ 
	\(z'=2yz+a_0\) \\
		\quad with \(f'' = 0\)
\end{tabular} & Painlev\'e II
  \\ \hline
\hyperref[sec:IX.5]{\(\text{IX.B}(5)\)} &
\begin{tabular}{l} 
	\(y'=-y^2+z+3q\) \\ 
	\(z'=4yz-9q'\) \\
		\quad with \((q''-6q^2)''=0\)
\end{tabular}  & Painlev\'e I
  \\ \hline
\hyperref[sec:11]{\(\text{XI}\)} &
\begin{tabular}{l} 
	\(y'= y(y-z)+ay\) \\ 
	\(z'= z(z-y)-az\)
\end{tabular} & ---
  \\ 
\hline 
\hyperref[sec:12]{\(\text{XII}\)} &
\begin{tabular}{l} 
	\(y'=y(y-2z)-2fy+b_0\) \\ 
	\(z'=z(z-2y)+2fz-a_0\) \\
	\quad with \(f''=0\) 
\end{tabular} & Painlev\'e IV
 \\ \hline
\hyperref[sec:13]{\(\text{XIII}\)} &
\begin{tabular}{l} 
	\(y'=\frac{1}{2}y(2z-y)+2py\) \\ 
	\(z'=\frac{1}{2}z(3y-2z)-4pz+2p^2-2p'+f\) \\
	\quad with \(f''=0\), \((p''-2p^3-fp)'=0\)
\end{tabular} & Painlev\'e II
  \\ \hline
\hyperref[sec:14]{\(\text{XIV}\)} &
\begin{tabular}{l} 
	\(y'=y(2z-y)+3py+r\) \\ 
	\(z'=z(y-z)-2pz\) \\
	 \quad with  \(r'=pr\), and, for \(q=\frac{1}{12}(p'+p^2-r)\),  \((q''-6q^2)''=0\)   
\end{tabular} & Painlev\'e I
   \\ \hline
\end{tabular} 
\caption{Representatives of the quadratic systems free of movable critical points (Theorem~\ref{thm:main}). Here, \(n\in \ZZ\), \(a\), \(b\), etc. are functions of the independent variable, \(a_0\) and \(b_0\) are constants.}\label{table:main} 
\end{table}

The labeling of the systems in Table~\ref{table:main} parallels Bureau's in \cite{bureau-9} and~\cite{bureau-10}. Systems \(\text{V}_0\), \(\text{VII}_0\), \(\text{IX.A}_0(n)\) and \(\text{IX.B}_0(n)\) in Table~\ref{table:main} are absent from Bureau's classification. They are simple equations, which are very easily integrated, but their absence challenges Bureau's claim of comprehensiveness in the presentation of Table~III in~\cite[\S 28]{bureau-10}: ``Il nous para\^\i t utile de r\'eunir dans le tableau suivant les repr\'esentants de \emph{tous} les s.d. [syst\`emes différentiels] polynomiaux stables de degré deux'' (we thought it useful to gather representatives of \emph{all} stable polynomial differential systems of degree two in the following table). We give slightly simpler formulas than Bureau's for systems~I, III and~IV, and correct an error in the analysis of system IX.A(3). Our normal forms for systems XII, XIII and XIV differ from Bureau's; they permit to relate them in a simpler way to relevant second-order equations in one variable, and to express in a more concise manner the conditions for the absence of movable critical points. We take into account the symmetries of systems XII and XIII, which allows us to avoid the repetition of some calculations (or, through it, to double or triple check their correctness). Apart from this, the results summarized in Table~\ref{table:main} are due to Bureau.
 
Some systems in  Table~\ref{table:main} are  birationally equivalent to the first-order system in two variables associated a Painlev\'e equation, and, through these identifications, some pairs of systems belonging to different rows in Table~\ref{table:main} can be shown to be  birationally equivalent. We have the following result, obtained  in the course of establishing Theorem~\ref{thm:main}.

\begin{proposition}\label{prop:bir_pain} For the systems in Table~\ref{table:main} that are solved by  Painlevé transcendents, a bimeromorphic map \(\phi:\CC^2\times\CC\dashrightarrow \CC^2\times\CC\), of the form \((\zeta,t)\mapsto (\xi(\zeta,t),t)\), with \(\xi\) birational as a function of \(\zeta\) for most values of \(t\), transforms them into the first-order system in two variables associated to the corresponding Painlev\'e equation. 
\begin{description}
\item[Painlev\'e I's birational class] The Painlev\'e I equation and some of its degenerations are given by the second-order equation
\begin{equation}\label{eq:PI-2nd}
u''=6u^2+f, \; \text{with } f''=0.
\end{equation}
Systems~V, IX.B(2),  IX.B(5) and XIV are birationally equivalent to the first-order system in two variables associated to it (and, in particular, are birationally equivalent to one another). More precisely, we have the following.
\begin{description}
\item[V] System~V is exactly the first-order system in two variables associated to~(\ref{eq:PI-2nd}), for the same parameter \(f\).
\item[IX.B(2)] From a solution \((y,z)\) to the system IX.B(2) with parameter~\(q\), \[u={\textstyle\frac{1}{6}}z+q\] 
is a solution to equation (\ref{eq:PI-2nd}) with parameter \(f=q''-6q^2\). (Since \(u\) is polynomially expressed in terms of \(y\) and \(z\), and the system solved by \(y\) and \(z\) is polynomial in \(y\) and \(z\), we have that \(u'\) may be polynomially expressed in terms of \(y\) and \(z\)). Reciprocally, for the equation (\ref{eq:PI-2nd}) with parameter \(f\), and a solution \(u\) to it, for \(q\) defined by \(q''=6q^2+f\),
\[(y,z)=\left(\frac{u'-q'}{u-q},6(u-q)\right)\]
is a solution to system IX.B(2) with parameter \(q\). These two transformations are inverses of one another.
\item[IX.B(5)] From a solution \((y,z)\) to the system IX.B(5) with parameter~\(q\),
\[u=y^2-{\textstyle\frac{1}{3}}z-2q,\]
is a solution to equation (\ref{eq:PI-2nd}) with parameter \(f=q''-6q^2\). Reciprocally, from a solution \(u\) to equation~(\ref{eq:PI-2nd}) with parameter \(f\),  for \(q\) defined by \(f=q''-6q^2\),
\[(y,z)=\left(-\frac{1}{2}\left(\frac{u'-q'}{u-q}\right),\frac{3}{4}\left(\frac{u'-q'}{u-q}\right)^2-3u-6q\right)\]
is a solution to system IX.B(5) with this parameter \(q\). These two transformations are inverses of one another.
\item[XIV] Consider a solution \((y,z) \) to system XIV with parameters \((r,p)\),  and let \(q=\frac{1}{12}(p'+p^2-r)\). Then,
\[u={\textstyle \frac{1}{6}}zy+{\textstyle \frac{1}{6}}r+q\]
is a solution to equation (\ref{eq:PI-2nd}) with parameter \(f=q''-6q^2\). Reciprocally, from a solution \(u\) to equation~(\ref{eq:PI-2nd}) with parameter \(f\), for \(q\) defined by \(q''-6q^2=f\), and \(p\) and \(r\) satisfying the conditions of system XIV with respect to \(q\),
\[(y,z)=\left(\frac{6u-6q-r}{\frac{u'-q'}{u-q}-p},\frac{u'-q'}{u-q}-p\right)\]
is a solution to system XIV with parameters \((r,p)\). These two transformations are inverses of one another.
\end{description}
\item[Painlev\'e II's birational class] The Painlev\'e~II equation and some of its degenerations are given by the second-order equation
\begin{equation}\label{eq:PII-2nd}
u''=2u^3+fu+\alpha, \; \text{with } \alpha\in \CC, f''=0.
\end{equation} 
Systems IX.B(3) and XIII  are birationally equivalent to the first order system in two variables associated to it.
\begin{description}
\item[IX.B(3)] From a solution \((y,z)\) to system IX.B(3) with parameters \((f,a_0)\),  \(y\) is a solution to~(\ref{eq:PII-2nd}) for the same \(f\) and \(\alpha=a_0-\frac{1}{2}f'\). Reciprocally, from  a solution \(u\) to the equation (\ref{eq:PII-2nd}) with parameters \((f,\alpha)\), 
\[(y,z)=\left(u,u'+u^2+{\textstyle\frac{1}{2}}f\right)\]
is a solution to the system IX.B(3) whose parameters are the same \(f\) and \(a_0=\alpha+\frac{1}{2}f'\). These two nonlinear polynomial transformations are inverses of one another.
\item[XIII] For  a solution \((y,z)\)  of system XIII with parameters \((f,p)\), 
\[u={\textstyle \frac{1}{2}}y-p\] 
is a solution to equation (\ref{eq:PII-2nd}) for the same \(f\) and
\begin{equation}\label{eq:rel:XII-PII}\alpha=-(p''-2p^3-fp).\end{equation} 
Reciprocally, given \(f\) with \(f''=0\) and \(\alpha\in \CC\), from a solution \(u\) to the associated system~(\ref{eq:PII-2nd}), for \(p\) defined by (\ref{eq:rel:XII-PII}),
\[(y,z)=\left(2(u+p),u-p+\frac{u'+p'}{u+p}\right)\]
is a solution to system XIII with parameters \((f,p)\). These two transformations are inverses of one another.
\end{description}
\item[Painlev\'e IV's birational class] The Painlev\'e IV equation and some of its degenerations are given, in Gambier's normalization \cite{gambier-acta}, by 
\begin{equation}\label{eq:PIV-2nd}
u''=\frac{(u')^2}{2u}+\frac{3}{2}u^3+4fu^2+2(f^2-\alpha)u+\frac{\beta}{u}, \; \text{with }  \alpha,\beta\in\CC, f''=0.
\end{equation}
System XII is birationally equivalent to the first-order system in two variables coming from it. If \((y,z)\) is a solution to system XII,  \(z\) 
solves Eq.~(\ref{eq:PIV-2nd}) with the same \(f\) and  \begin{equation}\label{eq:condbir-XII-P4}\alpha={\textstyle \frac{1}{2}}a_0+b_0-f', \;\beta=-{\textstyle\frac{1}{2}}a_0^2.\end{equation}  Reciprocally, from a solution \(u\) to equation (\ref{eq:PIV-2nd}) with parameters \((f,\alpha,\beta)\), for  constants  \(a_0\) and \(b_0\) satisfying~(\ref{eq:condbir-XII-P4}), 
\[(y,z)=\left(-\frac{u'}{2u}+\frac{u}{2}+f-\frac{a_0}{2u},u\right)\]
is a solution to system XII with these parameters.  
\end{description}

\end{proposition}

This implies that these systems are more than just free of movable critical points:
\begin{corollary} The solutions to the systems V, IX.B(2), IX.B(3), X.B(5), XII, XIII and XIV of Table~\ref{table:main} are single-valued meromorphic functions defined in all of~\(\CC \). 
\end{corollary}

Systems IX.B(3) and XIV precede Bureau's investigations, and already appear in Malmquist's formulation of Painlev\'e's equations II and IV as Hamiltonian systems \cite[\S 33]{malmquist-17}. The first-order systems directly associated to these equations are not Hamiltonian, but the  birational models given by the quadratic systems IX.B(3) and XIV are. This is presented in Table~\ref{tab:hamiltonian}, along with the Hamiltonian of the first Painlev\'e equation, directly related to system~V.

\begin{table}
	\begin{tabular}{lll}
Equation & Hamiltonian & system  	\\ \hline
Painlev\'e I & \(\frac{1}{2}z^2-2y^3-fy\) & V  \\
Painlev\'e II & \(-y^2z+\frac{1}{2}z^2-\frac{1}{2}fz-a_0y\)  &  IX.B(3)\\
Painlev\'e IV & \(y^2z-yz^2-2fyz+a_0y+b_0z\) & XII \\ \hline 
	\end{tabular}

\caption{Malmquist's Hamiltonian formulation for the Painlev\'e I, II and IV equations. Here,  \(a_0, b_0\in\CC\),  \(f''=0\).}\label{tab:hamiltonian}
\end{table}

The less satisfactory part of the classification concerns system VIII, the only one for which the necessary and sufficient conditions for the absence of movable critical points do not explicitly appear, neither in Bureau's work, nor in Table~\ref{table:main}. These systems form an infinite family, indexed by a positive integer \(n\), and the complexity of a certain necessary condition for the absence of movable critical points grows with~\(n\); it seems difficult to give an explicit expression for this condition for every~\(n\). At the same time, the equations are, conceptually, very simple, as they factor through Riccati equations.  In any case, there remains the question of whether this necessary condition is also a sufficient one. We will show that this is indeed the case, by studying the birational geometry of these systems, in the spirit of Okamoto's work on the spaces of initial conditions of the Painlev\'e equations~\cite{okamoto}, which is, on its turn, motivated by Painlev\'e's 17th of his \emph{Le\c cons de Stockholm}~\cite{painleve-lecons}. This will be done in Section~\ref{sec:VIII}.

\begin{remark}
The equations in Table~\ref{table:main} are not in their most reduced form. For a given system, there may be different parameter sets giving equivalent equations. We have favored the simplicity of the forms, at the cost of some minor redundancies. 
\end{remark}

\begin{remark}\label{rmk:noX}
No system X appears in Table~\ref{table:main}. A system X does appear in Bureau's Table~III in \cite{bureau-10}, but there are several issues with it.  First, in Bureau's rather structured text,  the discussion of any system labeled X would be expected to appear in \cite{bureau-9} after \S 20 but before~\S 21, and is nowhere to be found. Second, the system appearing in Bureau's table is \(y'=yz+a\), \(z'=(n-1)yz+b\), with \(n>1\); however, the parameter \(n\) in it is rendered irrelevant by scaling \(y\), which suggests that there is a misprint in the formula. Indeed, as per the caption for this system in Bureau's table, the quadratic part of the system should be the one in formula (19.3) in \cite{bureau-8}, which is  said to be equivalent to the quadratic system (19.2) there,  thus belonging to system~IX. We  suspect that our systems  \(\text{IX.A}_0(n)\) and \(\text{IX.B}_0(n)\) play the role conceived by Bureau for his system~X. In any case, the classification in Theorem~\ref{thm:main} has no need for this system.\end{remark}

This article can be read independently of Bureau's work, although we have included information aimed at facilitating the parallel reading of both.

\section{Criteria} 

We now describe some of the tools used by Bureau to establish the absence of movable critical points for the equations under consideration.

\subsection{Quadratic homogeneous equations and univalence} \label{sec:quad}

Consider a system of equations of the form (\ref{eq:int}), 
\[\zeta'=\sum_{i=0}^{n} P_{i}(\zeta,t),\]
with \(\zeta\in \CC^2\) 
and \(P_i\) a holomorphic function which, as a function of \(\zeta\), is a homogeneous polynomial of degree \(i\). For $\alpha\neq 0$, $\xi=\alpha z$, and $t=t_0+\alpha^{n-1} \tau$,
\[\frac{\dd\xi}{\dd\tau}=\sum_{i=0}^{n}\alpha^{n-i} P_{i}(\xi,t_0+\alpha \tau).\]
The limit as \(\alpha\to 0\) is the autonomous, homogeneous equation
\(\dd\xi/\dd\tau=P_n (\xi,t_0)\). If the original equation is free of movable critical points, the latter will have only single-valued solutions~\cite{painleve-bsmf} (see also~\cite{bureau-1}). Thus, in our setting, the first problem is to determine the equations of the form
\begin{equation}\label{eq:quadhom}
\begin{split}
y' & =  P(y,z),\\
z' & =  Q(y,z), 
\end{split} 
\end{equation} 
with \(P\) and \(Q\) quadratic homogeneous polynomials, having single-valued solutions. This has been done by Goffar-Lombet in~\cite[Ch.~I]{goffar-lombet}, by Bureau in~\cite{bureau-8}, and by Ghys and Rebelo in~\cite[Thm.~C]{ghys-rebelo-2} and~\cite{gr2-err} (and is likely to have been done elsewhere as well). The classification appears in Table~\ref{table:quad.hom}, where we have listed representatives for all the linear equivalence classes of univalent quadratic homogeneous equations of the form (\ref{eq:quadhom}); we follow Bureau's normalizations. Observe that these representatives have no continuous parameters in them.

\begin{table}
	\begin{tabular}{clc}
		System & reduced equation    & indices \\ \hline
		I &   \(\begin{array}{l}y'  = -y^2 \\ z'  = -yz \end{array}\) &   ---  \\ \hline
		V &  \(\begin{array}{l}y'  = 0  \\ z'  = 6y^2  \end{array}\) &   --- \\ \hline
		\begin{tabular}{c} III (\(n=-1\)) \\ IV (\(n=1\)) \\ IX (\(n\notin\{-1,0,1\}\)) \end{tabular}  &  \(\begin{array}{l}y'  =-y^2   \\ z'  =  (n-1)yz \end{array}\)    &  \(n\) \\  \hline
		VII &   \(\begin{array}{l}y'  = 0  \\ z'  =  yz \end{array}\)  & \(\infty\) \\ \hline
		\begin{tabular}{c}  II (\(n=0\)) \\ VIII (\(n>0\)) \end{tabular} &  \(\begin{array}{l}y'  = y[y-(n+2)z]  \\   z'  =  -z(ny+z) \end{array}\) &  \(1,n+1,-n-1\) \\ \hline
		VI &   \(\begin{array}{l}y'  = 0  \\ z'  = z(y+z) \end{array}\)  & \(1,\infty,\infty\) \\ \hline
		XI &  
		\(\begin{array}{l}y'  = y(y-z)  \\ z'  =  z(z-y) \end{array}\)
		&   \(2,2,\infty\)\\ \hline
		XII &   \(\begin{array}{l}y'  =  y(y-2z) \\ z'  =  z(z-2y) \end{array}\)    & \(3,3,3\) \\ \hline
		XIII & \(\begin{array}{l}y'= \frac{1}{2}y(2z-y)    \\ z'  =  \frac{1}{2}z(3y-2z) \end{array}\)  & \(2,4,4\) \\ \hline
		XIV &   \(\begin{array}{l}y'  =   y(2z-y)\\ z'  =  z(y-z) \end{array}\)   & \(2,3,6\) \\ 
		\hline
	\end{tabular}
\caption{Representatives of the univalent quadratic homogeneous systems, with the indices of their simple radial orbits  (we adhere to Bureau's normalization except in system XII, where our normal form differs from Bureau's by a change of sign in \(z\)). }\label{table:quad.hom}
\end{table}

An approach to this classification is the following. For a generic solution \((y,z)\) of~(\ref{eq:quadhom}), \(s=y/z\) is either constant, or satisfies an equation of the form 
\begin{equation}\label{schwarz-chris} 
-\frac{s''}{(s')^2}=F(s),
\end{equation} 
for some rational function \(F\) (a similar equation is satisfied by \(1/s\)). Briot and Bouquet gave necessary and sufficient conditions on \(F\) for the univalence of a function \(s\) satisfying such an equation~\cite{briot-bouquet}: at the finite values of \(s\), \(F\) has only simple poles, and their residues are either equal to \(-1\), or are of the form \(1/k-1\) for \(k\in\ZZ\), plus the analogous condition at \(s=\infty\) coming from the equation satisfied by \(1/s\) (see also Painlev\'e~\cite[\S 10]{painleve-bsmf}; see~\cite[Section~2.2]{guillot-quadratic} for an interpretation in terms of affine structures). The \emph{index} of the pole of (\ref{schwarz-chris}) is the associated integer \(k\) (we will say that a simple pole with residue \(-1\) has index \(\infty\)). We refer the reader to the closing remarks in~\cite[\S 2.3]{conte:painleve_approach} for a discussion around this terminology.

A \emph{radial orbit} for a homogeneous equation of the form (\ref{eq:quadhom}) is a line through \((0,0)\) that is invariant by the equation (that is a union of solutions). If not every line through the origin is a radial orbit, the equation has three of them (counted with multiplicity). The poles of the associated equation (\ref{schwarz-chris}) come from radial orbits, although some radial orbits may fail to produce them. The \emph{index} of a radial orbit is the one of the corresponding pole in equation (\ref{schwarz-chris}); a simple radial orbit which produces no poles for (\ref{schwarz-chris}) is said to have index~\(1\). If an equation has three radial orbits, with indices~\(n_i\), 
\[\frac{1}{n_1}+\frac{1}{n_2}+\frac{1}{n_3}=1.\]
The classification of the quadratic homogeneous equations (\ref{eq:quadhom}) with three different radial orbits that have only univalent solutions amounts to solving this equation as a Diophantine one (admitting that \(n_i\) can take the value \(\infty\)). The last six lines of Table~\ref{table:quad.hom} are obtained in this way (the linear equivalence class of such an equation is determined by its indices, and the equations thus obtained happen to be univalent). The other lines of this table  follow from organizing quadratic homogeneous vector fields according to the multiplicities of their radial orbits, to the indices of the simple ones, and from imposing the conditions guaranteeing the simplicity of the poles of the associated equation~(\ref{schwarz-chris}).

\subsection{Canonical equations and Malmquist's conditions} \label{sec:malmquist}

Consider the system of equations: 
\begin{equation}\label{malmeq}
\begin{split}
y\frac{\dd x}{\dd t} & =  na(t)x+b(t)y+\sum_{i+j\geq 2} f_{ij}(t)x^iy^j,\\  
\frac{\dd y}{\dd t} & =   a(t)+\sum_{i+j\geq 1} g_{ij}(t) x^iy^j,
\end{split}  
\end{equation} 
with holomorphic right-hand side, with \(n\in\ZZ\), \(n\geq 1\), and with \(a\) a non-vanishing function. Following Bureau, we will call these equations \emph{canonical}. The integer \(n\) appearing in such an equation will (also) be called its \emph{index}. We have the following result, appearing in~\cite[\S 11]{malmquist-2eme} (see also \cite[\S 2, Thm.~IV]{bureau-m1}), which was probably well-known before Malmquist.

\begin{proposition}[Malmquist] \label{prop:malm} In the case \(n=1\), for the system (\ref{malmeq}) to be free of movable critical points, it is necessary that \(b= 0\). \end{proposition}
We refer the reader to Malmquist's  work for a complete proof. The idea is, first, to establish, via a theorem of Dulac, that  for every \(t_0\)  there exists an invariant surface in the extended phase space of the form \(t=F(x,y)\) with \(F(0,0)=t_0\), and  that, in restriction to this surface, \(x\), as a function of \(y\), satisfies an equation of Briot-Bouquet type: \(y\, \dd x/\dd y=x+\cdots\); the solution of the latter will have logarithmic terms unless \(b(t_0)=0\). Since this should happen for all values of  \(t\), we have that \(b\) vanishes identically.

If \(n>1\), we may, in the extended phase space, blow up the curve \((x,y)=(0,0)\):  for \[u=\frac{b}{(n-1)a},\]
and for \(z\) defined by \(x=(z-u)y\), the system becomes
\begin{equation*}
\begin{split}
y\frac{\dd z}{\dd t} & =   (n-1)a\xi+yh(z,y),\\
\frac{\dd y}{\dd t} & =  a +\sum_{i+j\geq 1} g_{ij} (z-u)^iy^{i+j},
\end{split} 
\end{equation*}
for
\[h(z,y)=u'+\sum_{i+j\geq 2} f_{ij} (z-u)^iy^{i+j-2}-\sum_{i+j\geq 1} g_{ij} (z-u)^{i+1}y^{i+j-1}.\]
The transformed system is still a canonical one, but, in it, the index \(n\) has decreased by one. By iterating this procedure, we arrive to the case \(n=1\), where the criterion of Proposition~\ref{prop:malm} may be applied.

\section{The classification} 
We follow Bureau's steps for the greater part of the classification. We start with a system of equations of the form
\begin{equation}\label{test}
\begin{split}
y' & =  P(y,z,t)y+Ay+Bz+a, \\
z' & =  Q(y,z,t)z+Cy+Dz+b,
\end{split} 
\end{equation}
with \(P\) and \(Q\) quadratic homogeneous polynomials in \(y\) and \(z\) (at least one of them non-trivial), whose coefficients, along with \(A\), \(B\), \(C\), \(D\), \(a\) and \(b\), are holomorphic functions of the independent variable defined in a domain \(U\subset \CC\), and we suppose the system to be free of movable critical points. By the results in Section~\ref{sec:quad}, for all \(t\) in \(U\), the linear equivalence class of the quadratic part of~(\ref{test}) is represented by one of the autonomous systems of Table~\ref{table:quad.hom}, and is thus, with the exception of a discrete set of values of \(t\), locally constant (this is the first source of exceptional values of \(t\) in Theorem~\ref{thm:main}). By restricting \(t\) to a neighborhood of such a value, up to a change of coordinates that is linear in \(y\) and \(z\), and that preserves the independent variable, we may suppose that the quadratic part in (\ref{test}) is independent of \(t\), and that it is one of the quadratic systems of Table~\ref{table:quad.hom}. 

By applying suitable transformations that are affine on \(y\) and \(z\), possibly involving a change of independent variable, we will simplify the system without modifying its quadratic part. We will consider changes of coordinates of the forms
\begin{alignat*}{3}  
y  & = \lambda Y+h, \quad & z    &=\mu Z+k,   & t  &=\varphi(\tau);  \tag{\(\text{T}_1\)} \\
y  & =\lambda Y+h, \quad  & z  & =\lambda Z+k, &  t &  =\varphi(\tau);  \tag{\(\text{T}_2\)} \\
y  & = \lambda Y+h,\quad  & z   & =\mu Z+\nu Y+k, \quad & t  &=\varphi(\tau), \tag{\(\text{T}_3\)}
\end{alignat*}
where   \(\lambda\), \(\mu\), \(\nu \), \(h\) and \(k\) are functions of \(t\), chosen to solve a particular system of algebro-differential equations in a neighborhood of some \(t_0\in U\) (the need, in Theorem~\ref{thm:main}, to restrict the domains, is also a consequence of the  local character of this analysis).

For systems other than system~I,
when extending \(\CC^2\times U\) into \(\PP^2\times U\) by 
\begin{equation}\label{eq:intop2}
(y,z,t)\mapsto (y:z:1,t),
\end{equation}
the system will have poles along the line at infinity for all values of \(t\), with singularities at the directions of the radial orbits of the quadratic part. At a singularity associated to a radial orbit with positive index \(n\), we will obtain a canonical equation~(\ref{malmeq}) with index~\(n\), and the results in Section~\ref{sec:malmquist} will impose a necessary condition for the absence of movable critical points, of algebro-differential nature, which may be explicitly calculated after performing \(n-1\) blow-ups. We will meet some changes of coordinates that may be framed within this setting; for instance, the coordinates \((1/y,z)\) amount to considering, first, a change of affine chart in \(\PP^2\), \((u,v)=(1/y,z/y)\), followed by a blow-up, \((u,v)\mapsto (u,v/u)=(1/y,z)\).

In some cases, we will, after having  calculated all or some of these conditions, relate the system to a second-order equation in one variable, and join the Painlev\'e-Gambier classification at an advanced stage, either to obtain  conditions for the absence of movable critical points, or to integrate the equations. In this we differ from Bureau, who opts for calculating as many necessary conditions as possible directly from the system.

We will thus obtain systems with few parameters, which, in most cases, we will be able either to integrate, or to relate to a time-tested equation already known to be free of movable critical points.

We now begin the classification.  Systems III and IV will be studied together   with system~IX. As discussed in Remark~\ref{rmk:noX}, there is no system~X.

\subsection{System I}\label{sec:1}   Studied by Bureau in~\cite[\S 9]{bureau-9}. It represents the systems of the form 
\[
\begin{split}
y' & =  L(y,z,t)y+Ay+Bz+a, \\
z' & =  L(y,z,t)z+Cy+Dz+b,
\end{split}
\]
with \(L\) a linear form in \(y\) and  \(z\). These have long been considered as natural generalizations of the Riccati equation (see~\cite[\S 31]{malmquist-17}, \cite[Ch.~XXI]{goursat-CAMII}): when extending the phase space by (\ref{eq:intop2}), the system extends holomorphically. Their solutions are thus of the form
\[(y,z)=\left(\frac{\alpha y_0+\beta z_0+\gamma}{\eta y_0+\theta z_0+\iota},\frac{\delta y_0+\epsilon z_0+\zeta}{\eta y_0+\theta z_0+\iota}\right),\]
with \(\alpha\), \(\beta\), etc. functions of  \(t\), and these  equations are hence  naturally free of movable critical points throughout the domain where they are defined. A linear transformation redresses \(L(y,z,t)\) into  \(-y\), and, after this, a suitable transformation   \(\text{T}_3\)   brings the  system into
\[
\begin{split}
y' & =  -y^2+Bz, \\
z' & =   -yz  
\end{split} 
\]
(a form simpler than Bureau's, who uses the more restrictive transformation \(\text{T}_2\)). An explicit expression for its integration follows from the fact that  \((1/z)''=B\).

\subsection{System II}\label{sec:2}   Studied by Bureau in~\cite[\S 10]{bureau-9}. By means of a transformation  \(\text{T}_2\), the system may be cast into Bureau's form
\[
\begin{split}
y' & =  y(y-2z)+Ay+a, \\
z' & =   -z^2+Cy. 
\end{split}
\]
For \(y=1/u\) and \(z=v/u\) we have the canonical equation of index~\(1\)
\[
\begin{split}
u' & =  -1-Au+2v-au^2, \\
uv' & =  -v+Cu-Auv+v^2-au^2v. 
\end{split}
\]
From Proposition~\ref{prop:malm},  \(C=0\) is a necessary condition for the absence of movable critical points. We will  henceforth assume it. At the point \((u,v)=(0,1)\), the system is a canonical equation of index~\(1\), from which we obtain the necessary condition  \(A=0\). With it, \(u=y-z\) solves the  Riccati equation with holomorphic coefficients \(u'=u^2+a\), which is free of movable critical points. Thus, the system solved  by \(u\) and \(z\) (and hence the one solved by \(y\) and \(z\)) is free of movable critical points.

\subsection{System V}\label{sec:V}  Studied by Bureau in~\cite[\S 13]{bureau-9}. The general form of the system is
\begin{equation} \label{V.base}
 \begin{split}
y' & = Ay+Bz+a,\\
z' & = 6y^2+Cy+Dz+b.
\end{split}
 \end{equation}
Let us not spare details in order to discuss Bureau's omission in this case. The most general affine transformation preserving the systems~(\ref{V.base}) is of type  \(\text{T}_3\). By choosing it in the form 
\[y=\lambda Y+h, \quad z=\mu Z+\lambda\nu Y+k,  \quad \frac{\dd t }{\dd\tau}=\frac{\mu}{\lambda^2},\]
with \(\lambda\), \(\mu\), \(\nu\), \(h\) and \(k\) functions of \(t\),
 we obtain 
\begin{multline*} 
\frac{\dd Y}{\dd \tau} =\frac{\mu}{\lambda^2} \left(A+\nu B-\frac{\lambda'}{\lambda}\right)Y+\frac{\mu^2}{\lambda^3}BZ+\frac{\mu}{\lambda^3}(a+hA+kB-h'),\\
\frac{\dd Z}{\dd \tau} = 
6Y^2+
\frac{1}{\lambda}(C-\nu A-\nu^2 B+\nu D-\nu'+12h)Y+
\frac{\mu}{\lambda^2}\left(D-\nu B-\frac{\mu'}{\mu}\right)Z \\ +\frac{1}{\lambda^2}(b+hC+kD+6h^2-k'-\nu[a+hA+kB-h']). \end{multline*} 
(\(f'\) denotes the derivative of \(f\) with respect to \(t\).) This shows that, in~(\ref{V.base}), the condition \(B=0\) is, with respect to affine changes of coordinates, an intrinsic one. We leave out of our analysis the values of \(t\) where \(B\) has an isolated zero, and restrict   to domains where either \(B\) vanishes identically, or where it does not vanish at all (this is another source of special values of the independent variable in the statement of Theorem~\ref{thm:main}; conditions of the sort will be implicitly assumed throughout the rest of the article).

If \(B\neq 0\), the system may be brought  into Bureau's form
\begin{equation}\label{pain1-2d}
\begin{split}
y' & = z,\\
z' & = 6y^2+f, 
\end{split}
 \end{equation}
equivalent to the second-order equation \begin{equation}\label{P1}
y''=6y^2+f.
\end{equation}
The latter belongs to the class  leading to the first Painlev\'e equation,   studied in detail by Painlev\'e himself in~\cite[\S 15]{painleve-bsmf} (see also~\cite[\S 14.311, p.~328--329]{ince}). For it, Painlev\'e found that the necessary and sufficient condition for the absence of movable critical points   is \(f''=0\). If \(f'\neq 0\), the equation is, up to translations in the independent variable and scalings, the equation defining the first Painlev\'e transcendent \(y''=6y^2+t\).  If \(f\) is the constant \(f_0\), Eq.~(\ref{P1}) has the first integral \((y')^2=4y^3+2f_0y+K\), and is solved by elliptic functions and their rational degenerations. In this case, \(4y^3+2f_0y-z^2\) is a first integral for (\ref{pain1-2d}); its generic level surfaces are elliptic fibrations over the domain of the independent variable.

Bureau obtains the necessary condition \(f''=0\) through a canonical equation obtained from Eq.~(\ref{pain1-2d}). This works both ways: Bureau's analysis implies that, in Eq.~(\ref{P1}), \(f''=0\) is a necessary condition for the absence of movable critical points.  We will reobtain this condition in sections~\ref{sec:system-IX} and~\ref{sec:14}.
 
The case \(B=0\) is overlooked by Bureau. In it, in the original system~(\ref{V.base}), \(y\) solves  a linear equation with holomorphic coefficients, and so does \(z\); the system is free of movable critical points. A suitable transformation \(\text{T}_3\) reduces the system to  system \(\text{V}_0\) in Table~\ref{table:main}.

\subsection{System VI} \label{sec:6} Studied by Bureau in~\cite[\S 14]{bureau-9}.
It is a limit of system VIII as \(n\to\infty\). The system can be brought into Bureau's form
\begin{equation*}
\begin{split}
y' & = Bz+a,\\
z' & = z(z+y)+b,
\end{split}
\end{equation*}
by a transformation of type  \(\text{T}_2\). In the coordinates \(u=-1/z\) and \(v=y/z\), we find the canonical equation of index~\(1\),
\begin{equation*}
\begin{split}
u' & = 1+v+bu^2, \\
uv' & = v+Bu-au^2+v^2+bu^2v,
\end{split}
\end{equation*}
from which we obtain that \(B=0\) is a necessary condition for the absence of movable critical points. The condition is also sufficient: with it, \(z\) solves a Riccati equation with holomorphic coefficients.

\subsection{System VII} \label{sec:7} Studied by Bureau in~\cite[\S 16]{bureau-9}. It is a limit of system IX as \(n\to\infty\). Its general form is
\begin{equation}\label{VII.base}
 \begin{split}	
y' & =Ay+Bz+a,\\
z' & = yz+Cy+Dz+b.
\end{split}
\end{equation}
We consider changes of coordinates of the form \(\text{T}_1\), for which the condition \(B=0\) is an intrinsic one. 

If \(B\neq 0\), the system may be brought into Bureau's form
\[
\begin{split}
y' & = z+a, \\ 
z' & = yz+b.
\end{split} \]
In the coordinates \((u,v)\) defined by \(y=-1/u\), \(z=(v+\frac{1}{2})/u^2\), 
we find the canonical system of index~\(2\)
\[ \begin{split}
u' & = \textstyle\frac{1}{2}+v+au^2, \\ 
uv' & = v+au^2+2v^2+bu^3+2au^2v,
\end{split} \]
from which we obtain that \(a=0\) is a necessary condition for the absence of movable critical points. With it,		
\(y''=yy'+b\), and thus \(y\) solves the Riccati equation with holomorphic coefficients \(y'=\frac{1}{2}y^2+\int b\): the system is free of movable critical points.
		
If \(B=0\), in (\ref{VII.base}), \(y\) solves a linear equation with holomorphic coefficients, and so does \(z\) afterwards, and the system is free of movable critical points. The system can be redressed, via a transformation \(\text{T}_1\), into the system \(\text{VII}_0\) of Table~\ref{table:main}.

\subsection{System VIII}\label{sec:VIII}
Studied by Bureau in~\cite[\S 17]{bureau-9}. It has a parameter \(n\in\ZZ\), \(n\geq 1\). By means of a transformation \(\text{T}_2\), the system  may be brought into Bureau's form
\begin{equation}\label{VII-base}
\begin{split}
y' & = y[y-(n+2)z]+Ay+a, \\
z' & =-z(ny+z)+b.
\end{split}
\end{equation}
In the coordinates \((u,v)\) defined by
\(y=1/u\), \(z=(v+1)/u\), we have the canonical system of index~\(1\)
\begin{equation*}
\begin{split}
u' & = (n+1)-Au+(n+2)v-au^2,  \\
uv' & = (n+1)v-Au-Auv+(n+1)v^2+(b-a)u^2-au^2v.
\end{split}
\end{equation*}
A first necessary condition for the absence of movable critical points is thus \(A=0\), which we assume from now on. For \(y=-1/u\), we have the canonical equation of index~\(n\) 
\begin{equation}\label{VII-prepared}
\begin{split}
u' &  = 1+au^2+(n+2)uz,\\
uz' & = nz+bu-uz^2,
\end{split}
\end{equation}
and, through it, we find a necessary condition for each value of~\(n\). For the smallest ones, these are
\begin{alignat*}{2}
& \text{for } n=1, \quad && b=0; \\
& \text{for } n=2,\quad && b'=0; \\
& \text{for } n=3, && b''+ab-3b^2=0;\\  & \text{for } n=4, && b'''+4(a-4b)b'+2ba'=0;
\end{alignat*}
it seems difficult to give an explicit form of this condition for all values of \(n\).

In order to integrate the system (\ref{VII-base}) with \(A=0\),
Bureau transforms it into a particular case of the second-order equation in one variable studied by Gambier in \cite[Ch.~VII, Eq.~5, p.~53]{gambier-acta},  which he revisits in \cite[\S 49]{bureau-m1} (when \(a\neq 0\), he passes through the second-order equation satisfied by \(1/y\); when \(n\neq 2\) and \(b\neq 0\), through the one satisfied by \(z\); he deals separately with the case \(n=2\)). Gambier's equation has a rather convoluted scheme of integration. In it, the necessary and sufficient conditions for the absence of movable critical points are easy to establish for the equations corresponding to low values of \(n\); for higher values, either the absence of movable critical points can be established directly, through the vanishing of a particular explicit expression, or the equation may be transformed into one of the same family, in which \(n\) is replaced by \(n-2\). We refer to Gambier's work for details. 

The eventual movable critical points of the system are easy to locate. For \(z\) and \(w=y-z\), we have the system
\begin{equation*} 
\begin{split}
w' &  = w^2+a-b,\\
z' & = -(n+1)z^2-nwz+b.
\end{split}
\end{equation*}
In it, \(w\) solves a Riccati equation with holomorphic coefficients, and has thus meromorphic solutions, with movable poles; \(z\) solves a Riccati equation with meromorphic coefficients, in which the movable critical points, if any, come from the movable poles of~\(w\).

With the condition  \(A=0\), the necessary condition for the absence of movable critical points coming from the canonical equation (\ref{VII-prepared}) turns out to be a sufficient one.  We will establish this geometrically, without  integrating the equations, in   connection with the ideas in~\cite[Sect.~3.4]{ghys-rebelo-2}.  Embed \(\CC^2\) into \(\PP^1\times \PP^1\) via \((w,z)\to(w:1,z:1)\). Let \(u=1/w\) and \(v=1/z\), so that \((w:1,1:v)\), \((1:u,z:1)\) and \((1:u,1:v)\) give the three other charts of \(\PP^1\times \PP^1\). In these,
\begin{equation*}
\begin{split}  
w' & = w^2+a-b, \\
v' & = n+1+nvw-bv^2; 
\end{split} 
\end{equation*}
\begin{equation} \label{p1xp1-2}
\begin{split}
u' & =-1+(b-a)u^2, \\
uz' & =-nz+bu-(n+1)uz^2;
\end{split}
\end{equation}
\begin{equation} \label{p1xp1-3}
\begin{split}
u' & = -1+(b-a)u^2, \\
uv' & =  nv+(n+1)u-buv^2.
\end{split}
\end{equation}
The system has poles along \(u=0\), and nowhere else. In the chart (\ref{p1xp1-2}), the equation has a singularity with index \(n\) along  \(p:(u,z)=(0,0)\), and, in the chart~(\ref{p1xp1-3}), a singularity with index \(-n\) at   \(q:(u,v)=(0,0)\).

\emph{Elementary} transformations play an important role in the birational theory both of ruled surfaces and of the Riccati equation (see~\cite[Ch.~6, Sect.~7]{bhpv} and~\cite[Ch.~4, \S 1]{brunella};   they are called ``flips of the fiber'' in the latter). We will  perform   parametric (non-autonomous) versions of such transformations.

Start with the rational fibration \(\Pi:\PP^1\times \PP^1 \times U \to \PP^1\times U \), \((1:u,z:1,t)\mapsto (1:u,t)\). Let \(F\) be the fiber of \(\Pi\) above \(\{u=0\}\times U\). For the intersection with a slice \(\Sigma: t=t_0\), \(F\cap\Sigma\) is, within \(\Sigma\), a rational curve of vanishing self-intersection. The blowing up of \(p\) produces a surface \(E\), and transforms \(F\) into a surface \(F'\). Within the transform of \(\Sigma\), both \(E\cap\Sigma\) and \(F'\cap\Sigma\) are rational curves of self-intersection~\(-1\). By the parametric version of Grauert's criterion due to Nakano and Fujiki~\cite{nakano, fujiki-nakano}, \(F'\) may be contracted along the fibers of \(t\), producing a curve \(q'\) transverse to \(t\). This transforms \(E\) into a surface \(E'\), and, within the new transform of \(\Sigma\), \(E'\cap \Sigma\) is a curve of vanishing self-intersection. The total space \(\PP^1\times\PP^1\times U\) has been modified into a threefold \(M\), and the transform of \(\Pi\) subsists as a locally trivial rational bundle \(\Pi':M\to \PP^1\times U\). The fiber \(F\) has been replaced by \(E'\), but all the other fibers have gone by unchanged. The transformed equation has poles along \(E'\). There is a curve \(p'\) on \(E'\), along which we have a canonical equation with index \(n-1\); at the curve \(q'\), we find a canonical equation with index \(1-n\). This procedure can be iterated: when applied to a point of index~\(1\), the conditions of Proposition~\ref{prop:malm} are equivalent to the fact that the poles in the equation disappear after blowing up the point, thus removing all poles from the equation.

After iterating the procedure, \(\PP^1\times\PP^1\times U\) has been birationally transformed into a threefold where the equation extends holomorphically (where it no longer has poles), such that the fibers of the projection onto \(t\) compactify the original fiber \(\CC^2\times \{t\}\). This implies that the equation is free of movable critical points (see~\cite[\S 1]{okamoto}, where Okamoto resorts to Ehresmann's theorem, or Painlev\'e's 17th of his \emph{Le\c cons de Stockholm}~\cite{painleve-lecons}).

For example, when \(b=0\), starting from the 
chart \((1:u,z:1)\), after the
composition of \(n\) elementary transformations given by 
\(z=su^n\), the system transforms into 
\begin{equation*}
\begin{split}
u' & =-1-au^2, \\
s' & =nasu-(n+1)u^ns^2,
\end{split}
\end{equation*}
one which no longer has poles. (In particular, for all values of \(n\), the condition \(b=0\) guarantees the absence of movable critical points.)

\subsection{System IX (plus systems III and IV)} \label{sec:system-IX}
These were studied by Bureau in \S 11, \S 12 and \S 18 of \cite{bureau-9}. Their general form is
\begin{equation}\label{IX.base}
 \begin{split}
y' & = -y^2+Ay+Bz+a, \\
z' & = (n-1)yz+Cy+Dz+b,
\end{split}
\end{equation}
with \(n\in\ZZ\), \(n\neq 0\). System III corresponds to \(n=-1\), system IV to \(n=1\), and system X to the remaining cases. The latter is split in two: system IX.A corresponds to \(n<-1\), and system IX.B to \(n>1\). \emph{Warning: in system IX.A, Bureau's normalization differs from the one here used by the sign of~\(n\).} 

Let us be overly explicit in the calculations in this part, in which Bureau's analysis presents some oversights. Under a transformation of type \(\text{T}_1\), of the form 
\[y  = \lambda Y+h, \quad z =\mu Z+k, \quad  \frac{\dd t}{\dd \tau}=\frac{1}{\lambda},\] 
with \(\lambda\), \(\mu\), \(h\) and \(k\) functions of \(t\),
 we have
\begin{multline*} 
	\frac{\dd Y}{\dd \tau}  =   -Y^2+\frac{1}{\lambda}\left( A-\frac{\lambda'}{\lambda}-2h\right)Y+\frac{\mu}{\lambda^2}BZ+\frac{1}{\lambda^2}(a+hA+kB-h^2-h'),\\
	\frac{\dd Z}{\dd \tau}   =  (n-1)YZ+\frac{1}{\mu}( C+(n-1)k)Y+  \frac{1}{\lambda} \left(D-\frac{\mu'}{\mu}+(n-1)h\right)Z+ \\ +   \frac{1}{\lambda \mu}(b+(n-1)hk+h C+ kD-   k')
\end{multline*}
(\(f'\) denotes the derivative of \(f\) with respect to \(t\)).

\begin{remark}\label{IX.rmk} Under such changes of coordinates, we have that, in (\ref{IX.base}), 
\begin{itemize}
\item the condition \(B= 0\) is an intrinsic one; and that
\item  if \(n=-3\) and \(B\neq 0\),  so is the  condition
\(2A-D-B'/B=0\).
\end{itemize}
\end{remark}

Our analysis will be split into the following cases: 

\subsubsection{Case \(n=1\) (system IV)}\label{sec:IV}
The system can be brought, by a transformation \(\text{T}_1\) as above, into the form
\begin{equation} \label{IV}
 \begin{split}
y' & = -y^2+Bz, \\
z' & = Cy
\end{split}
 \end{equation}
(a form simpler than Bureau's, who uses the more restrictive transformation \(\text{T}_2\)). In the coordinates \(u=1/y\), \(v=z/y\), it reads
\begin{equation*}
 \begin{split}
u' & = 1-Bvu, \\
uv' & = v+Cu-Buv^2,
\end{split}
 \end{equation*}
a canonical equation of index~\(1\), from which we obtain the necessary condition \(C=0\). With it, system (\ref{IV}) reduces to a Riccati equation in \(y\) with holomorphic coefficients, and is thus free of movable critical points.

We will henceforth suppose that \(n\neq 1\).

\subsubsection{Case \(B\neq 0\), \(n\neq-3\)} Up to a suitable transformation as above, we may suppose that the system has Bureau's form
\begin{equation} \label{IX.BnZ.qhs}
\begin{split}
y' & =-y^2+z+a,  \\
z' & =(n-1)yz+b. 
\end{split} 
\end{equation}
When \(y\) and \(z\) are respectively given the weights \(1\) and \(2\), its quasihomogeneous component of highest degree is
\begin{equation} \label{IX.qh}
\begin{split}
y' & =-y^2+z,  \\
z' & =(n-1)yz.
\end{split} 
\end{equation}

In this way, (\ref{IX.BnZ.qhs}) can be seen as a lower-order perturbation of~(\ref{IX.qh}), and the whole procedure described for quadratic equations can be adapted to the study of such systems. For instance, if (\ref{IX.BnZ.qhs}) is free of movable critical points, (\ref{IX.qh}) has only single-valued solutions. For \(v=z/y^2\), we have, from (\ref{IX.qh}), the Briot-Bouquet equation
\begin{equation}\label{eq:inv_aff}
-\frac{v''}{(v')^2}=-\frac{3v-n}{v(2v-n-1)},\end{equation}
discussed in Section~\ref{sec:quad}. We must impose the condition \(n\neq -1\) in order for the pole at \(v=0\) to be simple. The expression (\ref{eq:inv_aff}) has a simple pole at  \(v=v_1\) for 
\[v_1={\textstyle\frac{1}{2}}(n+1),\]
with residue 
\[\frac{n-1}{2(n+1)}-1.\]
If \(v\) is univalent, the latter is either equal to \(-1\) or has the form \(1/k-1\) for \(k\in \ZZ\), this is, \(2(n+1)/(n-1)=2+4/(n-1)\) is an integer (the index). The admissible cases are \(n\in\{2,3, 5\}\). The indices of the poles of (\ref{eq:inv_aff}) appear in Table~\ref{table:ind-qh21}.

\begin{table}
\begin{tabular}{c|ccc}
\(n\)   & 0 & \(v_1\) & \(\infty\) \\ \hline
\(2\) & \(3\) & \(6\)& \(2\)  \\ 
\(3\) & \(4\)& \(4\) & \(2\)  \\
\(5\) & \(6\) &  \(3\) & \(2\) \\ \hline  
\end{tabular}
\caption{Indices of the Briot-Bouquet equation associated to  equation (\ref{IX.qh}).}\label{table:ind-qh21}
\end{table}

Considering, with Bureau, the coordinates \(u=1/y\), \(v=z/y^2\), we have
\begin{equation} \label{IX:bur_cond}
\begin{split}
u' & = 1-v-au^2,\\
uv'  & =  (n+1)v-2v^2+bu^3-2au^2b.
\end{split}
 \end{equation}
We have two singular points along \(u=0\): \(v=0\) and \(v=v_1\). At each one of them, we have a canonical equation whose index appears in Table~\ref{table:ind-qh21}. This will give, in each of the admissible cases, two necessary conditions for the absence of movable critical points.

\paragraph{Subcase \(n=2\)} \label{sec:IX.2} 
The necessary condition coming from~(\ref{IX:bur_cond}) along \((u,v)=(0,0)\) is \(b=0\), which we assume from now on. Set \(a=12 q\). For \(u=\frac{1}{6}z+q\), we obtain \(u''-6u^2=q''-6q^2\), the equation appearing in Section~\ref{sec:V} as Eq.~(\ref{P1}) for \(f=q''-6q^2\), related to the first Painlev\'e equation. As we there discussed, the necessary and sufficient condition for the absence of movable critical points is \(f''=0\); this is also the necessary condition arising from Eq.~(\ref{IX:bur_cond}) at \((u,v)=(0,v_1)\). 

The system is birationally equivalent to the first-order system in two variables associated to the first Painlev\'e equation (cf. Proposition~\ref{prop:bir_pain}).

\paragraph{Subcase \(n=3\)}\label{sec:IX(3)}
We have, for \(f=-2a\) and \(g=b+a'\), 
\begin{equation} \label{P2}
y''=2y^3+fy+g.
\end{equation}
We find in~\cite[\S 14.315, p.~333]{ince} that the conditions for this equation to be free of movable critical points are \(f''=0\) and \(g'=0\) (\(a''=0\) and \(b'=0\)). These are exactly the two conditions arising from Eq.~(\ref{IX:bur_cond}) (they will also be  discussed in Section~\ref{sec:13}). Under these, (\ref{P2}) is the general form of the equation leading to the Painlev\'e~II transcendent~\cite[\S 11]{painleve-acta}. If \(f'\neq 0\), up to a translation the independent variable and suitable scalings, Eq.~(\ref{P2}) gives one of the standard forms of the Painlev\'e~II equation,
\[y''=2y^3+ty+g_0.\]
When \(f'=0\), (\ref{P2}) is integrated by elliptic functions and their rational degenerations: if \(f\) and \(g\) are, respectively, the constants \(f_0\) and \(g_0\), we have the first integral \((y')^2=y^4+f_0y^2+2g_0y+K\) (see also~\cite[\S 14.315, p.~333]{ince}). When substituting for \(y'\) from (\ref{IX.BnZ.qhs}), we obtain that  \(2y^2z-z^2+2by-2az-a^2\) is a first integral for (\ref{IX.BnZ.qhs}). Its generic level surfaces are elliptic fibrations over the domain of the independent variable.

The system has the Hamiltonian
\[-y^2z+{\textstyle\frac{1}{2}z^2+az-by},\]
equivalent to Malmquist's Hamiltonian for the Painlev\'e II equation \cite[\S 33]{malmquist-17} (we recover Malmquist's normalization in the coordinates \((\upsilon,\zeta)=(z+a,y)\)). The system is polynomially equivalent to the first-order system  associated to the Painlev\'e~II equation (cf. Proposition~\ref{prop:bir_pain}).

\paragraph{Subcase \(n=5\)} \label{sec:IX.5} The necessary condition arising from Eq.~({\ref{IX:bur_cond}) at \((u,v)=(0,v_1)\) is \(b=-3a'\). Set \(a=3q\) and \(b=-9q'\). We have that \(w=-2y\) solves the equation
\begin{equation}\label{eq:XIV-2nd} w''=-ww'+w^3-12qw+12q'.\end{equation}
This equation is integrated in \cite[p.~23, Eq.~6]{painleve-acta}: if \(u\) is a solution to 
\begin{equation}\label{altpain1}
u''-6u^2=q''-6q^2,
\end{equation}
and is different from \(q\),
\begin{equation}\label{eq:smartchange}
w=\frac{u'-q'}{u-q}
\end{equation}
is a solution to~(\ref{eq:XIV-2nd}). When considered as first-order equations in two variables, equations~(\ref{eq:XIV-2nd}) and (\ref{altpain1}) are birationally equivalent, since, from (\ref{eq:smartchange}), \(u\) and \(u'\) can be expressed rationally in terms of \(w\) and \(w'\):
\begin{equation*}
\begin{split}
	u & =  \textstyle{\frac{1}{6}} w^2-q+\textstyle{\frac{1}{6}} w', \\ 
	u' & = \textstyle{\frac{1}{6}} w^3-2qw+\textstyle{\frac{1}{6}}ww'+q',
\end{split} 
\end{equation*}
and vice versa.

As discussed in Section~\ref{sec:V}, the necessary and sufficient condition for (\ref{altpain1}) to be free of movable critical points is 
\begin{equation} \label{eq:enq}
(q''-6q^2)''=0,
\end{equation}
in which (\ref{altpain1}) becomes the first Painlev\'e equation or one of its degenerations. Thus, condition (\ref{eq:enq})  is a necessary and sufficient condition for (\ref{eq:XIV-2nd}) to be free of movable critical points as well (compare with \cite[\S 14.314, pp.~331--333]{ince}). Condition (\ref{eq:enq}) is also the one arising from the canonical equation (\ref{IX:bur_cond}) at \((u,v)=(0,0)\). 

The system is birationally equivalent to the first-order one  associated to the first Painlev\'e equation (cf. Proposition~\ref{prop:bir_pain}).

\subsubsection{Case \(B\neq 0\), \(n=-3\)} \label{sec:IX.-3} In this case, owing to the condition given in the second item of Remark~\ref{IX.rmk}, the system (\ref{IX.base}) cannot always, despite Bureau's claim, be brought into the form (\ref{IX.BnZ.qhs}). It can nevertheless be put into the form
\begin{equation*}
 \begin{split}
y' & =-y^2+z,\\
z' & =-4yz+Dz+b.
\end{split} 
\end{equation*}
For \(u=1/y\), \(v=z/y^2+1\), 
\begin{equation*} 
 \begin{split}
u'&=2-w, \\
uw'&=2v-Du+Duv-2v^2+bu^3,
\end{split}
 \end{equation*}
and we thus find the necessary condition \(D=0\) (Bureau's analysis of this case fails to see that there are instances in which the system is not free of movable critical points). We have
\[y''+6yy'+4y^3-b=0.\]
For \(y=\frac{1}{2}w'/w\), this equation reduces to the third-order linear equation \(w'''=2bw\), and is thus free of movable critical points.

\subsubsection{Case \(B= 0\)}\label{sec:IX_0} This case was overlooked by Bureau. Through a transformation \(\text{T}_1\), the system may be brought into the form
\begin{equation*}
\begin{split}
y' & =-y^2,\\
z' & =(n-1)yz+Cy.
\end{split}
 \end{equation*}
After solving for \(y\), the system reduces to a linear inhomogeneous equation in \(z\) with meromorphic coefficients. Its solutions may have movable critical points, arising from the poles of \(y\). Its non-constant solutions are given by \(y(t)=1/(t-t_0)\) and
\[z(t)=(t-t_0)^{n-1}\int^t\frac{C(s)\,\dd s}{(s-t_0)^n}.\]
Thus, for \(n<0\), the system is always free of movable critical points. For \(n> 1\), the condition for the absence of movable critical points is \(\dd^{n-1}C/\dd t^{n-1}=0\). This gives the normal forms for systems \(\text{IX.A}_0(n)\) and \(\text{IX.B}_0(n)\) in Table~\ref{table:main}.

System III corresponds to the case \(n=-1\). In it, by using a transformation \(\text{T}_1\) instead of one of type \(\text{T}_2\), we have obtained a normal form simpler than Bureau's. The one appearing in Table~\ref{table:main} is obtained by changing the sign of~\(y\).

\subsection{System XI} \label{sec:11} Studied by Bureau in~\cite[\S 21]{bureau-9}.
A transformation \(\text{T}_2\) brings the system into the form
\begin{equation}\label{XI.base}
 \begin{split}
y' & = y(y-z)+qy+a, \\
z' & = z(z-y)-qz+b.
\end{split}
 \end{equation}
(This normal form is simpler than Bureau's, and exhibits the natural symmetries of the system.) For \(u=-1/y\), we have the canonical equation of index~\(1\)
\begin{equation*}
\begin{split}
u' & = 1-qu+uz+au^2, \\
uz' & = z+bu-quz+uz^2,
\end{split}
 \end{equation*}
from which we obtain the necessary condition \(b=0\). By the symmetry of the family~(\ref{XI.base}) consisting in exchanging \(y\) and \(z\), we have that \(a=0\) is a necessary condition as well. With these conditions, the system has the first integral \(K=yz\), and reduces to the family of Riccati equations with holomorphic coefficients \(y'=y^2+qy-K\). It is free of movable critical points.

\subsection{System XII}\label{sec:12}  Studied by Bureau in \cite[\S 22]{bureau-9}. By means of a transformation \(\text{T}_2\), the system may be brought into the form
\begin{equation}\label{XII.form}
 \begin{split}
		y'  & =  y(y-2z)-2fy+b,\\
		z' & = z(z-2y)+2fz-a. 
	\end{split}
 \end{equation} 
This form is slightly different from Bureau's, already at the level of the quadratic homogeneous part, which differs from the one in \cite[\S 22]{bureau-9} by a change of sign in \(z\). System (\ref{XII.form}) is Hamiltonian: it may be written as
\(y'=\partial \Phi/\partial z\), \(z'=-\partial \Phi/\partial y\),  
for 
\begin{equation}\label{for:hamiltonian_PIV}
\Phi=y^2z-yz^2-2fyz+ay+bz.
\end{equation}

For \(u=-1/y\),  we have
\begin{equation*}
\begin{split}
		u' & = 1+2fu+bu^2+2uz, \\
		uz' & = 2z-au+2fuz+uz^2,
	\end{split}
\end{equation*}
a canonical equation of index~\(2\), from which we obtain the necessary condition for the absence of movable critical points  \(a'=0\). The symmetry of the family~(\ref{XII.form}) consisting of interchanging \(y\) and \(z\)  and acts upon its coefficients by
\[(f,a,b)\mapsto (-f,-b,-a);\]
from this, we obtain \(b=0\) as a second  necessary condition for the absence of movable critical points. Under the change of coordinates
\((y,z)\mapsto(-y,z-y+2f)\),   system (\ref{XII.form}) is mapped into  one of the same family, with the transformation upon the coefficients  given by
\[(f,a,b)\mapsto (-f,a+b-2f',-b).\]
This is involutive. Through this, we find that \(f''=0\) is a third necessary condition for the absence of movable critical points. With these three conditions, \(z\) is a solution to  Eq.~(\ref{eq:PIV-2nd}) for \(\alpha=\frac{1}{2}a+b-f'\), \(\beta=-\frac{1}{2}a^2\).

Eq.~(\ref{eq:PIV-2nd}) is free of movable critical points: it is the general form of the equation leading to the Painlev\'e IV equation~\cite[Ch.~1, \S 9]{gambier-acta}. When \(f'\neq 0\), up to translations in the independent variable and scalings, we may suppose that \(f(t)=t\), recovering the Painlev\'e IV equation in its full generality. The cases where \(f'=0\) are integrated by elliptic functions: if \(f\) is the constant \(f_0\), Eq.~(\ref{eq:PIV-2nd}) has the first integral \(K\) defined by
\[(y')^2=y^4+4f_0y^3+4(f_0^2-\alpha)y^2+Ky-2\beta\]
(see also~\cite[\S 14.331, case 6\(^\circ\), p.~339]{ince}). In this case, in which \(f\) is constant, the system (\ref{XII.form}) is an autonomous one, and has in (\ref{for:hamiltonian_PIV}) a first integral; its generic level curves are elliptic.

System (\ref{XII.form}), with the integrability conditions, is birationally equivalent to the first-order system in two variables associated to the Painlev\'e IV equation~(\ref{eq:PIV-2nd}); see Proposition~\ref{prop:bir_pain}. The Hamiltonian (\ref{for:hamiltonian_PIV}) is the one associated by Malmquist to the fourth Painlev\'e equation \cite[\S 33]{malmquist-17}.

The above involutions  generate a group of order six, isomorphic to the group \(S_3\) of permutations in three symbols.

\subsection{System XIII}\label{sec:13} 
Studied by Bureau in~\cite[\S 23--24]{bureau-10}. Up to a transformation \(\text{T}_2\), we may suppose that the system has the form
\begin{equation}\label{XIII.form}
\begin{split}
y'  & =  \textstyle{\frac{1}{2}}y(2z-y)+2py+a,\\
z' & = \textstyle{\frac{1}{2}}z(3y-2z)-4pz+b.
	\end{split}
\end{equation}
(It is different from Bureau's.) For \(u = 1/z\), we have the canonical equation of index~\(1\)
\begin{equation*} 
\begin{split}
		u'  & = \textstyle 1+4pu-bu^2-\frac{3}{2}uy, \\
		uy' & = \textstyle y+au+2puy-\frac{1}{2}uy^2.
	\end{split}
\end{equation*}
Hence, \(a=0\) is a necessary condition for the absence of movable critical points. With it,
\(\frac{1}{2}y-p\)
solves Eq.~(\ref{P2}) for
\begin{equation}\label{xii.fandg}
\begin{split}
		f & = 2p'-2p^2+b,\\
		g & = -p''+2pp'+bp.
	\end{split}
\end{equation}
Following the discussion in Section~\ref{sec:IX(3)}, the necessary and sufficient conditions for Eq.~(\ref{P2}) to be free of movable critical points are \(f''=0\) and \(g'=0\), in which case we obtain the Painlev\'e~II equation and some of its degenerations. From~(\ref{xii.fandg}),
\begin{equation*}
\begin{split}
	b & = f-2p'+2p^2, \\
		p'' & = 2p^3+fp-g.
	\end{split}
\end{equation*}
In particular, with the above conditions, \(p\) is a solution to another instance of Eq.~(\ref{P2}), one with the same \(f\) but with an opposite sign for~\(g\).

The system is birationally equivalent to the first-order one associated to the Painlev\'e II equation (see Proposition~\ref{prop:bir_pain}).
 
The conditions \(f''=0\) and \(g'=0\) can be obtained through some canonical equations obtained from system (\ref{XIII.form}) with \(a=0\). For \(u=2/y\),  we have the canonical equation with index~\(3\) 
\begin{equation*}
 \begin{split}
		u' & = 1-2up-uz,\\
		uz' & = 3z+bu-uz^2-4puz,
	\end{split}
 \end{equation*}
from which we obtain the necessary condition for the absence of movable critical points \(f''+2g'=0\). 

The systems (\ref{XIII.form}) with \(a=0\) are preserved by the involution \((y,z)\mapsto (-y,z-y+4p)\), which acts upon their coefficients by \((p,b)\mapsto (-p,b+4p')\). The involution acts upon \(f\) and \(g\) by \((f,g)\mapsto (f,-g)\). Hence, we also have the necessary condition for the absence of movable critical points \(f''-2g'=0\). Together, these conditions give the previously discussed ones \(f''=0\) and \(g'=0\).

\subsection{System XIV}\label{sec:14} Studied by Bureau in~\cite[\S 25--27]{bureau-10}.
A transformation \(\text{T}_2\) brings the system to the form 
\begin{equation}\label{XIV.form} 
\begin{split}
		y'  & =  y(2z-y)+3py+a,\\
		z' & =  z(y-z)-2pz+b 
\end{split}
\end{equation}
(different from Bureau's). For \( (u,v)=(1/y,z/y)\), we have the canonical equation of index~\(2\)
\begin{equation} \label{eq:14.aux}
\begin{split}
		u'  & =  1-3pu-2v-au^2,\\
		uz' & =  2v+bu^2-5puv-3v^2-au^2v,
\end{split}
\end{equation}
from which we obtain the necessary condition \(b=0\). With it, for \(u=1/z\), we have the canonical system of index~\(2\)
\begin{equation*}
\begin{split}
		u'  & = 1+2pu-uv ,\\
		uv' & = 2v+au+3puv-uv^2 ,
\end{split}
\end{equation*}
and, from it, a second necessary condition, \(a'=pa\), follows. Assuming these two, we find that \(z+p\) is a solution to~(\ref{eq:XIV-2nd}) for \(q=\frac{1}{12}(p'+p^2-a)\). By eliminating \(a\) from \(q\) and \(q'\), we have that \(p\) is a particular solution to~(\ref{eq:XIV-2nd}). System~(\ref{XIV.form}) and the first-order system associated to~(\ref{eq:XIV-2nd}) are birationally equivalent. As discussed in Section~\ref{sec:IX.5}, the necessary and sufficient condition for Eq.~(\ref{eq:XIV-2nd}) to be free of movable critical points is \((q''-6q^2)''=0\). Hence, for system~(\ref{XIV.form}), assuming the first two necessary conditions, this is a necessary and sufficient condition as well. 

The system is birationally equivalent to the first-order one associated to the first Painlev\'e equation (cf. Proposition~\ref{prop:bir_pain}).

At \((u,v)=(0,2/3)\), Eq.~(\ref{eq:14.aux}) is a canonical system of index~\(6\), and, assuming the first two necessary conditions for the absence of movable critical points, the condition arising from this point is \((q''-6q^2)''=0\); we may thus obtain the necessity of this condition directly from system~(\ref{XIV.form}).

\subsection*{Acknowledgments} The author thanks the referees for their helpful comments.


\providecommand{\doi}[1]{\url{https://doi.org/#1}}

\end{document}